\documentclass[preprint,12pt]{elsarticle}

\usepackage{amssymb}
\usepackage{multirow}
\usepackage{booktabs}
\usepackage{amssymb}
\usepackage{mathrsfs}
\usepackage{longtable}
\usepackage{lscape}
\usepackage{tabularx}
\usepackage{epsfig}
\usepackage{colortbl}
\usepackage{subfigure}
\usepackage{amsmath}
\usepackage{inputenc}
\usepackage{float}
\usepackage{enumerate}
\usepackage{multirow}
\usepackage{rotating}
\usepackage{xcolor}
\usepackage{geometry}
\newtheorem{myDef}{Definition}
\journal{Annals of Operations Research}

\begin{document}

\begin{frontmatter}



\title{A DEMATEL-Based Completion Method for Incomplete Pairwise Comparison Matrix in AHP}


\author[1]{Xinyi Zhou}
\author[2]{Yong Hu}
\author[1,5,2,6]{Yong Deng\corref{cor} }
\cortext[cor]{Corresponding author: Yong Deng, School of Computer and Information Science, Southwest University, Chongqing, 400715, China.}
\ead{ydeng@swu.edu.cn, ydeng@xjtu.edu.cn, prof.deng@hotmail.com}
\author[3]{Felix T.S. Chan}
\author[4]{Alessio Ishizaka}

\address[1]{School of Computer and Information Science, Southwest University, Chongqing, 400715, China.}

\address[2]{Big Data Decision Institute, Jinan University, Tianhe, Guangzhou, 510632, China.}

\address[5]{Institute of Integrated Automation, School of Electronic and Information Engineering, Xi'an Jiaotong University, Xi’an, Shaanxi, 710049, China.}

\address[6]{School of Engineering, Vanderbilt University, Nashville, TN, 37235, USA.}

\address[3]{Department of Industrial and Systems Engineering, The Hong Kong Polytechnic University, Hung Hom, Hong Kong}

\address[4]{University of Portsmouth, Portsmouth Business School, Richmond Building, Portland Street, Portsmouth PO1 3DE, United Kingdom}

\begin{abstract}

Pairwise comparison matrix as a crucial component of AHP, presents the preference relations among alternatives. However, in many cases, the pairwise comparison matrix is difficult to complete, which obstructs the subsequent operations of the classical AHP. In this paper, based on DEMATEL which has ability to derive the total relation matrix from direct relation matrix, a new completion method for incomplete pairwise comparison matrix is proposed. The proposed method provides a new perspective to estimate the missing values with explicit physical meaning. Besides, the proposed method has low computational cost. This promising method has a wide application in multi-criteria decision-making.

\end{abstract}
\end{frontmatter}

\section{Introduction}
\label{Introduction}

Analytic hierarchy process (AHP) is a multi-criteria decision-making method that helps the decision-makers facing a complex problem with multiple conflicting and subjective criteria\cite{ishizaka2009analytic}.
AHP has been widely applied in various areas\cite{zhu2016hesitant,lin2015measuring,yepes2015cognitive,wang2011linear}.
Based on a hierarchical structure, the criteria weights can be derived from the pairwise comparison matrix.
Pairwise comparison matrix as a crucial component of AHP, is commonly utilized to estimate preference values of finite alternatives with respect to a given criterion. 

However, the pairwise comparison matrix is difficult to complete in the following cases, which 
obstructs the subsequent operations of the classical AHP.


\begin{itemize}
\item The experts lack the knowledge of one or more alternatives.
\item Partial data in pairwise comparison matrix has been lost.
\item Huge number of pairwise comparison requested: $(n^2-n)/2$ for n alternatives. For example, 8 alternatives and 6 criteria require 183 entries\cite{ishizaka2011review}.
\end{itemize}

Hitherto, it has attracted widespread attention over the issue of incomplete pairwise comparison matrix in AHP.
There are quantities of methods having been developed to address the incomplete comparison matrices. 
Some of them use the incomplete matrix to assign priority weight to each alternative without completion\cite{harker1987incomplete,carmone1997monte,hua2008ds,xu2013logarithmic,dopazo2011parametric}. Some methods are devoted to estimate the missing values in matrix\cite{fedrizzi2007incomplete,ergu2016estimating,bozoki2010optimal,chen2015bridging,zhang2012linear,hu2006backpropagation,alonso2008consistency,benitez2015consistent,gomez2010estimation,buyukozkan2012new,cabrerizo2010managing}.
While these methods have a good effect on allocating the priority weights to alternatives, the methods that focus on matrix completion are more ``well-rounded'' to restore the preference relations between each pair of alternatives, which are inaccessible before completion.

After reviewing the existing completion methods for incomplete pairwise comparison matrix in AHP, it can be found that the vast majority of these methods focus on the optimal consistency of pairwise comparison matrix\cite{fedrizzi2007incomplete,ergu2016estimating,bozoki2010optimal,chen2015bridging,zhang2012linear}. 
For example,
Fedrizzi and Giove\cite{fedrizzi2007incomplete} proposed a new method to complete the pairwise comparison matrix by minimizing its global inconsistency.
Ergu \emph{et al.}\cite{ergu2016estimating} extended the geometric mean induced bias matrix to estimate the missing judgments and improve the consistency ratios by the least absolute error method and the least square method.
Based on the graph theory, Bozo\'iki \emph{et al.}\cite{bozoki2010optimal} developed a completion method which focuses on the optimal consistency and estimates the missing values through iteration.
Besides,
Hu and Tsai\cite{hu2006backpropagation} proposed a well-known back propagation multi-layer perception to estimate the missing comparisons of incomplete pairwise comparison matrices in AHP.
Alonso \emph{et al.}\cite{alonso2008consistency} defined a way to calculate the consistency of incomplete pairwise comparison matrix and maintain its consistency level after completion.
Ben\'itez \emph{et al.}\cite{benitez2015consistent} used the Frobenius norm to obtain the most consistent/similar complete matrix with the initial incomplete matrix.
Gomez \emph{et al.}\cite{gomez2010estimation} proposed a completion method based on the neural network. These methods are useful but have some limitations:


\begin{itemize}
\item Most methods complete the pairwise comparison matrix with optimal consistency, which sometimes can not restore the incomplete matrix to its original state, especially when the initial matrix is not high consistent.
\item The methods are computationally expensive.
\end{itemize}

Taking into account these issues, a DEMATEL-based completion method for incomplete pairwise comparison matrix in AHP is proposed.
Based on DEMATEL method which has ability to derive the total relation matrix from direct relation matrix, the completion of incomplete pairwise comparison matrix can be divided into three steps. 

\begin{itemize}
\item[-] Firstly, turn the incomplete comparison matrix into direct relation matrix. 
\item[-] Secondly, DEMATEL method is adopted to obtain the total relation matrix. 
\item[-] Thirdly, transform the total relation matrix into pairwise comparison matrix with reciprocal preference relations. 
\end{itemize}

The proposed method provides a new perspective to estimate the missing values in pairwise comparison matrix with explicit physical meaning (i.e. by calculating the total relations among alternatives). Besides, the proposed method has low computational cost.

The rest of this paper is organized as follows. Section 2 introduces some preliminaries of this work. Section 3 illustrates the procedure of the proposed method and an application for ranking tennis players is presented. Section 4 discusses the superiority and efficiency of the proposed method and section 5 ends the paper with the conclusion.


\section{Preliminaries}
\label{Rreparatorywork}

\subsection{Analytic Hierarchy Process (AHP)}

AHP is a multi-criteria decision-making method developed by Saaty\cite{saaty1980analytic} and aims at quantifying relative weights for a given set of criteria on a ratio scale\cite{xu2015deviation}.
As a decision-making approach, it permits a hierarchical structure of criteria, which provides users with a better focus on specific criteria and sub-criteria when allocating the weights\cite{ishizaka2013multi}.
Besides, AHP has ability to judge the consistency of pairwise comparison matrix to show its potential conflicts in decision-making process.
AHP has been widely used in
supply chain management\cite{chan2007global,deng2014supplier,chan2006ahp}, 
human reliability analysis\cite{su2015dependence},
cloud computing\cite{liu2016decision},
energy planning\cite{ishizaka2016energy,onar2015multi,kahraman2009comparative},
site selection\cite{erdougan2016combined,lee2015integrated,pourahmad2015combination},
software assessment\cite{hu2013software,hu2015cost} and
performance evaluation\cite{zavadskas2015selecting,wu2010evaluating}.
AHP can be extended further by other theories and methods such as fuzzy number\cite{chan2008global,wang2008extent}, choquet integral\cite{corrente2015combining}, k-means algorithm\cite{lolli2014new}, VIKOR\cite{buyukozkan2015evaluation} and TOPSIS\cite{erdogan2014evaluating,kahraman2015fuzzy}.

Generally, the process of applying AHP is divided into three steps. 
Firstly, establish a hierarchical structure by recursively decomposing the decision problem.
Secondly, construct the pairwise comparison matrix to indicate the relative importance of alternatives.
A numerical rating including nine rank scales is suggested, as shown in Tab.\ref{Rating}. 
Thirdly, verify the consistency of pairwise comparison matrix and calculate the priority weights of alternatives. For the completeness of explanation, several basic concepts and formulas are introduced as follows.

\begin{table}[H]
\centering
\caption{Numerical rating in AHP}
\label{Rating}
 \footnotesize
\begin{tabular*}{0.65\textwidth}{@{\extracolsep{\fill}}ll}
 \toprule[1pt]
\quad Scale & Meaning \\ \hline
\quad 1 & Equal importance \\
\quad 3 & Moderate importance \\
\quad 5 & Strong importance \\
\quad 7 & Demonstrated importance \\
\quad 9 & Extreme importance \\
\quad 2, 4, 6, 8 & Intermediate values \\
\bottomrule[1pt]
\end{tabular*}
\end{table}

\begin{myDef}
(Pairwise comparison matrix)
Assume $\begin{Bmatrix} E_1, E_2, \cdots, E_n \end{Bmatrix}$ are $n$ alternatives available for decision-making, the pairwise comparison matrix is indicated as $M=\left ( m_{ij} \right )_{n \times n} \left (i,j=1,2,\cdots ,n\right )$, and satisfies:
\begin{equation}
\label{PCM}
m_{ij}=\left\{\begin{matrix}
\frac{1}{m{ji}} & i\neq j\\ 
1 & i=j
\end{matrix}\right.
\end{equation}
where $m_{ij}$ represents the relative importance of $E_{i}$ over $E_{j}$.
\end{myDef}

Consistency checking is introduced in AHP to verify the usability of constructed pairwise comparison matrix.

\begin{myDef}
(Consistency ratio)
For the pairwise comparison matrix $M_{n \times n}$, 
let $\lambda _{\max }$ denote the largest eigenvalue of $M$, consistency index ($CI$) is defined as
\begin{equation}
\label{CI}
CI = \frac{{\lambda _{\max }  - n}}{{n - 1}}
\end{equation}
Based on $CI$, consistency ratio (CR) is defined as 

\begin{equation}
\label{CR}
CR = \frac{{CI}}{{RI}}
\end{equation}
where $RI$ is the random consistency index related to the dimension of matrix, listed in Tab.\ref{RI}.

\begin{table}[H]
\centering
\caption{Random consistency index $RI$}
\label{RI}
 \footnotesize
\begin{tabular*}{1\textwidth}{@{\extracolsep{\fill}}cccccccccccc}
\toprule[1pt]
 \multicolumn{2}{c}{n} & 1 & 2 & 3 & 4 & 5 & 6 & 7 & 8 & 9  &10  \\ \hline
\multicolumn{2}{c}{$RI$} & 0 & 0 & 0.52 & 0.89 & 1.12 & 1.26 & 1.36 & 1.41 & 1.46 & 1.49  \\
\bottomrule[1pt]
\end{tabular*}
\end{table}

\end{myDef}

If $CR<0.1$, the constructed pairwise comparison matrix is considered acceptable and the priority weights of alternatives can be obtained through the eigenvector as Def.\ref{eigenvector}.
Otherwise, the comparison matrix needs to be reconstructed.

\begin{myDef}
\label{eigenvector}
(Eigenvector)
For the pairwise comparison matrix $M_{n \times n}$ with acceptable consistency, 
suppose $\vec{w} = ( {w_1 ,w_2 ,\cdots,w_n })^T $ is the eigenvector of $M$, whose $w_i(i=1, 2, \cdots, n)$  is indicated as the priority weight of the $i$-th alternative and calculated by
\begin{equation}
\label{weights of AHP}
M\vec{w} = \lambda _{\max }\vec{w}
\end{equation}
\end{myDef}

\subsection{DEMATEL method}
Decision-making trial and evaluation laboratory (DEMATEL) method was originally developed by Battelle Memorial Institute of Geneva Research Center, aimed at the fragmented and antagonistic phenomena of world societies and searched for integrated solution\cite{gabus1972world,fontela1976dematel}.
DEMATEL is built on the basis of graph theory, 
effectively assessing the total relations and constructing a map between system components in respect to their type and severity\cite{tzeng2007evaluating}.
Reviewing the former studies, DEMATEL has been successfully applied in many diverse areas such as 
service quality analysis\cite{shieh2010dematel,tseng2009causal},
stock selection\cite{shen2015combined} and
supply chain management\cite{ ISI:000368039100008,wu2015case,wu2015exploring}.
DEMATEL can be further extended by other theories and methods such as fuzzy numbers\cite{mentes2015fsa,tsai2015using} and ANP\cite{ ISI:000356990800033,tzeng2012combined, ISI:000262178000049,wu2008choosing}.

The procedure of DEMATEL is divided into five steps. 

\emph{\textbf{Step 1:} Define quality feature and establish measurement scale}

Quality feature is a set of influential characteristics that impact the sophisticated system, which can be determined by literature review, brain-storming, or expert evaluation. 
After defining the influential characteristics in researching system, establish the measurement scale for the causal relationship and pairwise comparison among influential characteristics.
Four levels 0,1,2,3 are suggested, respectively meaning ``no impact'', ``low impact'', ``high impact'' and ``extreme high impact''. In this step, factors and their direct relations are displayed by a weighted and directed graph.


\emph{\textbf{Step 2:} Extract the direct relation matrix of influential factors} 

In this step, transformation from the weighted and directed graph into direct relation matrix has been done.
For $n$ influential characteristics $F_1, F_2, \cdots, F_n$, direct relation matrix is denoted as $D=\left(d_{ij}\right)_{n \times n} \left (i,j=1,2,\cdots ,n\right )$, 
where $d_{ij}$ is the direct relation of $F_i$ over $F_j$ based on the measurement scale and satisfies $d_{ij}=0$ if $i=j$.

\emph{\textbf{Step 3:} Normalize the direct relation matrix}

The normalized direct relations of factors are a mapping from $d_{ij}$ to $\left [ 0,1 \right ]$, defined as below.

\begin{myDef}
\label{NorDRM}
(Normalized direct relation matrix)
For the framework of $n$ influential characteristics $\begin{Bmatrix} F_1, F_2, \cdots, F_n \end{Bmatrix}$, normalized matrix $N$ of direct relation matrix $D=\left(d_{ij}\right)_{n \times n} \left (i,j=1,2, \cdots ,n\right )$ is obtained by
\begin{equation}
\label{normMatrix}
N = \frac{D}{\max (\mathop { } \sum\limits_{j = 1}^n {{d_{ij}},} \mathop {} \sum\limits_{i = 1}^n {{d_{ij}}} )}
\end{equation}
\end{myDef}

Note that there must be at least one row of matrix $D$ satisfying $\sum_{j=1}^{n}d_{ij}<{\max (\mathop { } \sum\limits_{j = 1}^n {{d_{ij}},}\mathop {} \sum\limits_{i = 1}^n {{d_{ij}}} )} $\cite{lin2008causal}. Hence,
based on the studies of \cite{papoulis1985probability} and \cite{goodman1988introduction}, sub-stochastic matrix can be obtained through utilizing the normalized direct relation matrix $N$
and
absorbing state of Markov chain matrices, shown in (\ref{normMatrix2}).
\begin{equation}
\label{normMatrix2}
\mathop {\lim }\limits_{k \to \infty } N^k=O \quad and \quad \mathop {\lim }\limits_{k \to \infty } (I + N + N^2  +  \cdots  +N^k ) = N(I - N)^{ - 1}
\end{equation}
where $O$ is the null matrix and $I$ is the identity matrix.

\emph{\textbf{Step 4:} Calculate the total relation matrix}

Due to the characteristic of normalized direct relation matrix $N$, total relation matrix which contains direct and indirect relations among factors can be derived from (\ref{TRM}). 

\begin{myDef}
(Total relation matrix)
For the framework of $n$ influential characteristics $\begin{Bmatrix} F_1, F_2, \cdots, F_n \end{Bmatrix}$, assume $N=\left(n_{ij}\right)_{n \times n} \left (i,j=1,2,\cdots ,n\right )$ is the normalized direct relation matrix, total relation matrix $T$ is defined as
\begin{equation}
\label{TRM}
T =\mathop {\lim }\limits_{k \to \infty } (N + N^2  +  \cdots  +N^k )= N(I - N)^{ - 1}
\end{equation}
where $I$ is the identity matrix.
\end{myDef} 

\emph{\textbf{Step 5:} Construct the DEMATEL map}
 
\begin{figure}[H]
\begin{center}
\psfig{file=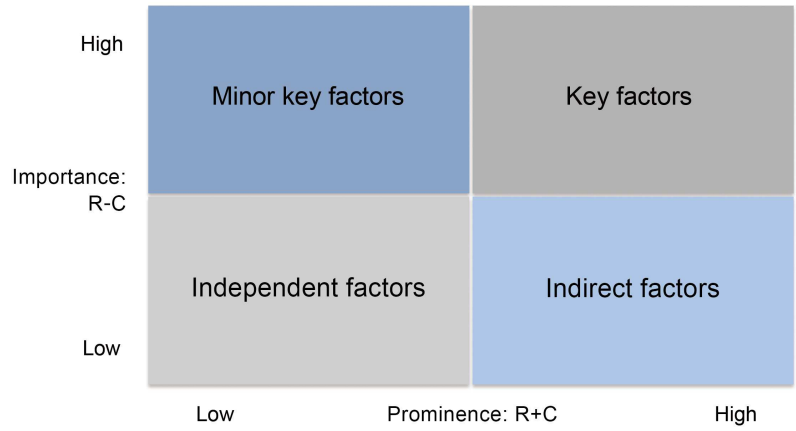,scale=0.56,angle=0}
\caption{DEMATEL map}
\label{DEMATELdiagram}
\end{center}
\end{figure}

Based on the sum of each row $R_i(i=1,2,\cdots,n)$ and column $C_i$ $(i=1,2,\cdots,n)$ of the total relation matrix $T_{n \times n}$, $(R_i+C_i)$ and $(R_i-C_i)$ 
can be obtained. 
$(R_i+C_i)$ is defined as the prominence, indicating the importance of the i-th influential factor.
$(R_i-C_i)$ classifies the i-th influential factor into the cause or effect category in researching system. If $(R_i-C_i)>0$, the factor is regarded as a cause factor. 
If $(R_i-C_i)<0$, the factor is an effect one. The DEMATEL map is shown as Fig.\ref{DEMATELdiagram}.


\begin{figure}[H]
\begin{center}
\psfig{file=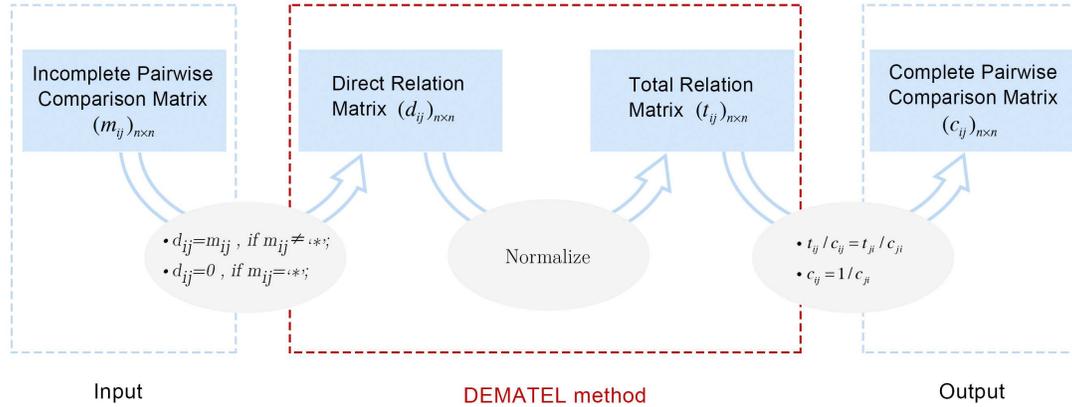,scale=0.51,angle=0}
\caption{Procedure of DEMATEL-based completion method}
\label{flowchart}
\end{center}
\end{figure}



\newpage
\section{DEMATEL-based completion method}
\label{ProposedMethod}
In this section, a DEMATEL-based completion method for incomplete pairwise comparison matrix in AHP is proposed. Moreover, an application for ranking top tennis players is presented.

\subsection{Procedure of DEMATEL-based completion method}


In this sub-section, the procedure of DEMATEL-based completion method is presented and illustrated by an example. Fig.\ref{flowchart} shows the process of the proposed method.

Assume an $4\times4$ incomplete pairwise comparison matrix $M$ and the unavailable/missing values are displayed by `*'.

\begin{displaymath}
M = \left[ {\begin{array}{*{20}c}
1 & * & 4 & 8\\ 
* & 1 & 2 & 4\\ 
1/4 & 1/2 & 1 & 2 \\
1/8 & 1/4 & 1/2 & 1 
\end{array}}\right]
\end{displaymath}

\emph{\textbf{Step 1:} From incomplete pairwise comparison matrix to direct relation matrix.}

Pairwise comparison matrix reflects the preference relations between each pair of factors,
where $m_{ij}(i,j=1,2,3,4)$ in $M$ indicates the relative importance of factor $i$ to factor $j$. Hence, the known values in $M$ can be put into the direct relation matrix $D$ without any transformation and the unavailable ones(`*') are substituted by 0. The direct relation matrix $D$ is as:

\begin{displaymath}
D = \left[ {\begin{array}{*{30}c}
1 & 0 & 4 & 8\\ 
0 & 1 & 2 & 4\\ 
1/4 & 1/2 & 1 & 2 \\
1/8 & 1/4 & 1/2 & 1 
\end{array}}\right]
\end{displaymath}

\emph{\textbf{Step 2:} From direct relation matrix to total relation matrix. }

First of all, normalize the direct relation matrix $D$ according to (\ref{normMatrix}). Then, based on the normalized matrix $N$ and (\ref{TRM}), total relation matrix $T$ can be obtained. The normalized direct relation matrix $N$ and total relation matrix $T$ are calculated as:
\begin{displaymath}
N = \left[ {\begin{array}{*{30}c}
0.0667 & 	0.0000 &	0.2667 &	0.5333 \\
0.0000 &	0.0667 &	0.1333 &	0.2667 \\
0.0167 &	0.0333 &	0.0667 &	0.1333 \\
0.0083 &	0.0167 &	0.0333 &	0.0667 
\end{array}}\right]
,
T = \left[ {\begin{array}{*{30}c}
0.0835 &	0.0241 &	0.3371 &	0.6742 \\
0.0060 &	0.0835 &	0.1685 &	0.3371 \\
0.0211 &	0.0421 &	0.0899 &	0.1798 \\
0.0105 &	0.0211 &	0.0449 &	0.0899 
\end{array}}\right]
\end{displaymath}

\emph{\textbf{Step 3:} From total relation matrix to pairwise comparison matrix.}

In fact, the total relation matrix has completed the unavailable values in initial incomplete matrix. 
Nevertheless, in pairwise comparison matrix, the multiplication of each pair of symmetric values along the diagonal line is required to equal 1, yet the values in total relation matrix are between 0 and 1.
Based on (\ref{transTotal2PCM}), it is feasible to accomplish this transformation from the total relation matrix $T=(t_{ij})_{n\times n}(i,j=1,2,\cdots,n)$ to pairwise comparison matrix $M_c=(c_{ij})_{n\times n} (i,j=1,2,\cdots,n)$.

\begin{equation}
\label{transTotal2PCM}
\left\{ {\begin{array}{*{20}l}
\frac{t_{ij}}{c_{ij}}=\frac{t_{ji}}{c_{ji}}  \\
c_{ij}=\frac{1}{c_{ji}} 
\end{array}} \right.
\end{equation}

Hence, the pairwise comparison matrix transformed from $T_{4 \times 4}$ is:
\begin{displaymath}
M_c = \left[ {\begin{array}{*{30}c}
1.0000 &	2.0000 &	4.0000 &	8.0000 \\
0.5000 &	1.0000 &	2.0000 &	4.0000 \\
0.2500 &	0.5000 &	1.0000 &	2.0000 \\
0.1250 &	0.2500 &	0.5000 &	1.0000 
\end{array}}\right]
\end{displaymath}

Finally, based on the proposed method, the completion for incomplete pairwise comparison matrix $M$ has been accomplished as:

\begin{displaymath}
M = \quad
\left[ {\begin{array}{*{20}c}
1 & * & 4 & 8\\ 
* & 1 & 2 & 4\\ 
1/4 & 1/2 & 1 & 2 \\
1/8 & 1/4 & 1/2 & 1 
\end{array}}\right]
\rightarrow 
\left[ {\begin{array}{*{20}c}
1 & \textbf{2} & 4 & 8 \\ 
\textbf{1/2} & 1 & 2 & 4\\ 
1/4 & 1/2 & 1 & 2 \\
1/8 & 1/4 & 1/2 & 1 
\end{array}}\right]
\end{displaymath}

\subsection{An application for ranking tennis players}


Professional tennis as a global popular sport owns large numbers of viewers around the world.
The professional tennis associations such as ATP have been collecting data about the tournaments and the players.
There is a free available database that collects the results of the top tennis players since 1973 (see at http://www.atpw-\\orldtour.com/Players/Head-To-Head.aspx).

\begin{figure}[htbp]
\begin{center}
\psfig{file=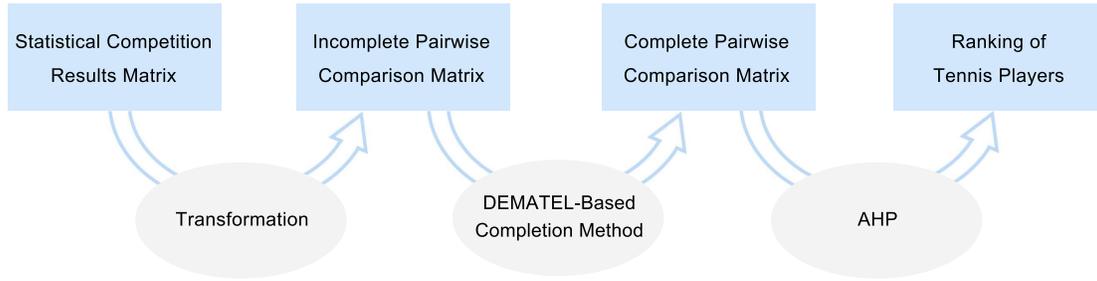,scale=0.16,angle=0}
\caption{Flowchart for ranking tennis players}
\label{rankFlowchart}
\end{center}
\end{figure}

In original statistics, 25 tennis players who have been the champion on ATP since 1973 are selected. The statistical competition results matrix is shown in Tab.\ref{statisticalMatrix}.
However, the statistical matrix is not complete since some of players did not have the chance to compete with each other or the partial data has been lost. Thus, DEMATEL-based completion method is adopted to rank the tennis players based on these statistics. The process is divided into three steps, which is shown in Fig.\ref{rankFlowchart}.




\begin{figure}[htbp]
\begin{center}
\psfig{file=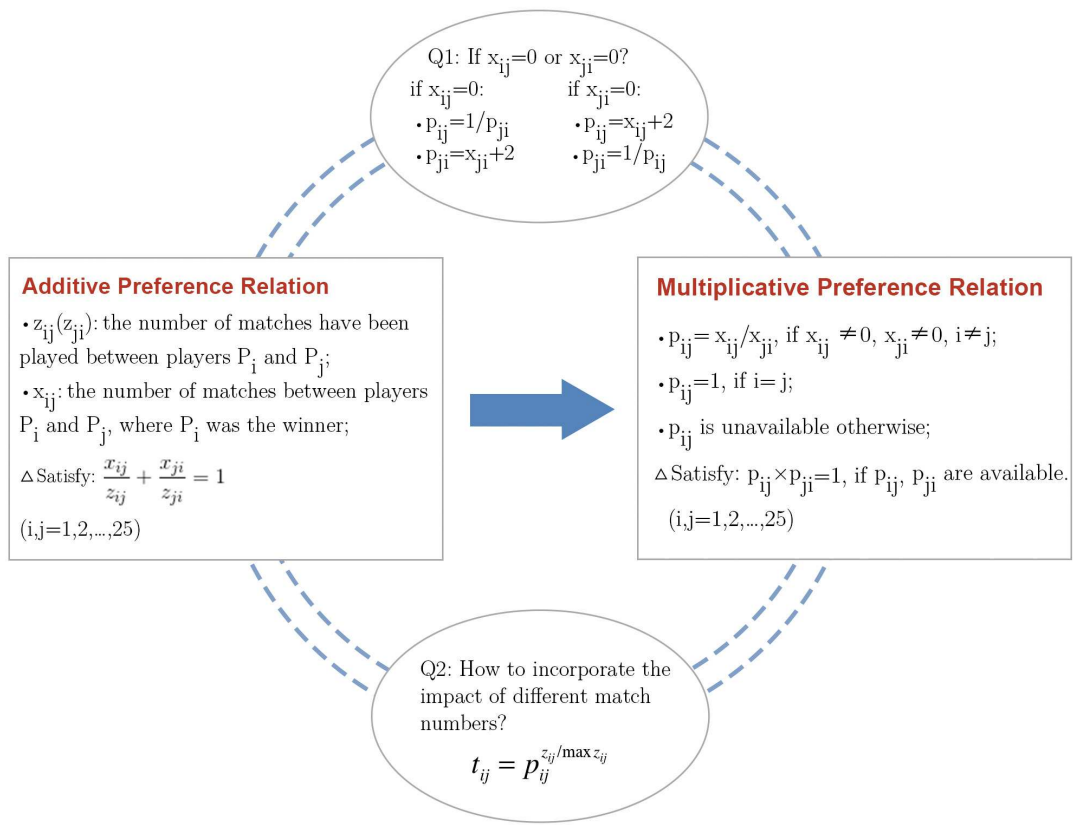,scale=0.68,angle=0}
\caption{Rule of transformation}
\label{transformation}
\end{center}
\end{figure}





\newgeometry{left=3cm,bottom=0.4cm,top=0.4cm}
\begin{sidewaystable}
\centering
\caption{Statistical competition results matrix: the wins/total number of matches}
\label{statisticalMatrix}
 \tiny
\begin{tabular*}{1\textwidth}{@{\extracolsep{\fill}}llllllllllllllllllllllllll}
\toprule[1pt]
           & Aga & Bec & Bor & Con & Cou & Djo & Edb & Fed & Fer & Hew & Kaf & Kue & Len & McE & Moy & Mus & Nad & Nas & New & Raf & Rio & Rod & Saf & Sam & Wil \\ \hline
Agassi     &        & 10/14  &      & 2/2     & 5/12    &          & 6/9    & 3/11    & 2/5     & 4/8    & 8/12       & 7/11    & 2/8   & 2/4     & 3/4  & 5/9    & 0/2   &         &          & 10/15  & 1/3  & 5/6     & 3/6   & 14/34   & 5/7      \\
Becker     & 4/14   &        &      & 6/6     & 6/7     &          & 25/35  &         &         & 1/1    & 4/6        &         & 10/21 & 8/10    & 2/4  & 2/3    &       & 1/1     &          & 2/3    & 3/5  &         & 0/1   & 7/19    & 7/10     \\
Borg       &        &        &      & 15/23   &         &          &        &         &         &        &            &         & 6/8   & 7/14    &      &        &       & 10/15   & 1/4      &        &      &         &       &         & 1/1      \\
Connors    & 0/2    & 0/6    & 8/23 &         & 0/3     &          & 6/12   &         &         &        &            &         & 13/34 & 14/34   &      &        &       & 12/27   & 2/4      &        &      &         &       & 0/2     & 0/5      \\
Courier    & 7/12   & 1/7    &      & 3/3     &         &          & 6/10   &         &         &        & 1/6        & 1/1     & 0/4   & 2/3     & 2/3  & 7/12   &       & 0/1     &          & 0/3    & 0/3  &         & 1/2   & 4/20    &          \\
Djokovic   &        &        &      &         &         &          &        & 15/31   & 2/3     & 6/7    &            &         &       &         & 2/4  &        & 17/39 &         &          &        &      & 4/9     & 0/2   &         &          \\
Edberg     & 3/9    & 10/35  &      & 6/12    & 4/10    &          &        &         &         &        & 1/3        &         & 14/27 & 6/13    & 1/1  & 10/10  &       &         &          & 3/3    & 1/1  &         &       & 6/14    & 9/20     \\
Federer    & 8/11   &        &      &         &         & 16/31    &        &         & 10/13   & 18/26  & 2/6        & 1/3     &       &         & 7/7  &        & 10/32 &         &          & 0/3    & 2/2  & 21/24   & 10/12 & 1/1     &          \\
Ferrero    & 3/5    &        &      &         &         & 1/3      &        & 3/13    &         & 4/10   & 1/3        & 3/5     &       &         & 8/14 &        & 2/9   &         &          & 2/3    & 3/4  & 0/5     & 6/12  &         &          \\
Hewitt     & 4/8    & 0/1    &      &         &         & 1/7      &        & 8/26    & 6/10    &        & 7/8        & 3/4     &       &         & 7/12 &        & 4/10  &         &          & 3/4    & 3/5  & 7/14    & 7/14  & 5/9     &          \\
Kafelnikov & 4/12   & 2/6    &      &         & 5/6     &          & 2/3    & 4/6     & 2/3     & 1/8    &            & 5/12    &       &         & 3/6  & 1/5    &       &         &          & 3/5    & 6/8  &         & 2/4   & 2/13    & 1/2      \\
Kuerten    & 4/11   &        &      &         & 0/1     &          &        & 2/3     & 2/5     & 1/4    & 7/12       &         &       &         & 4/7  & 3/3    &       &         &          & 4/8    & 2/4  & 1/2     & 4/7   & 1/3     &          \\
Lendl      & 6/8    & 11/21  & 2/8  & 21/34   & 4/4     &          & 13/27  &         &         &        &            &         &       & 21/36   &      & 4/5    &       & 1/1     &          & 0/1    &      &         &       & 3/8     & 15/22    \\
McEnroe    & 2/4    & 2/10   & 7/14 & 20/34   & 1/3     &          & 7/13   &         &         &        &            &         & 15/36 &         &      &        &       & 6/9     & 1/2      &        &      &         &       & 0/3     & 7/13     \\
Moya       & 1/4    & 2/4    &      &         & 1/3     & 2/4      & 0/1    & 0/7     & 6/14    & 5/12   & 3/6        & 3/7     &       &         &      & 4/8    & 2/8   &         &          & 3/4    & 2/7  & 1/5     & 4/7   & 1/4     &          \\
Muster     & 4/9    & 1/3    &      &         & 5/12    &          & 0/10   &         &         &        & 4/5        & 0/3     & 1/5   &         & 4/8  &        &       &         &          & 0/3    & 3/4  &         & 0/1   & 2/11    & 0/2      \\
Nad      & 2/2    &        &      &         &         & 22/39    &        & 22/32   & 7/9     & 6/10   &            &         &       &         & 6/8  &        &       &         &          &        &      & 7/10    & 2/2   &         &          \\
Nastase    &        & 0/1    & 5/15 & 15/27   & 1/1     &          &        &         &         &        &            &         & 0/1   & 3/9     &      &        &       &         & 4/5      &        &      &         &       &         & 0/1      \\
Newcombe   &        &        & 3/4  & 2/4     &         &          &        &         &         &        &            &         &       & 1/2     &      &        &       & 1/5     &          &        &      &         &       &         &          \\
Rafter     & 5/15   & 1/3    &      &         & 3/3     &          & 0/3    & 3/3     & 1/3     & 1/4    & 2/5        & 4/8     & 1/1   &         & 1/4  & 3/3    &       &         &          &        & 2/3  &         & 1/1   & 4/16    & 1/3      \\
Rios       & 2/3    & 2/5    &      &         & 3/3     &          & 0/1    & 0/2     & 1/4     & 2/5    & 2/8        & 2/4     &       &         & 5/7  & 1/4    &       &         &          & 1/3    &      & 0/2     & 1/4   & 0/2     &          \\
Roddick    & 1/6    &        &      &         &         & 5/9      &        & 3/24    & 5/5     & 7/14   &            & 1/2     &       &         & 4/5  &        & 3/10  &         &          &        & 2/2  &         & 4/7   & 2/3     &          \\
Safin      & 3/6    & 1/1    &      &         & 1/2     & 2/2      &        & 2/12    & 6/12    & 7/14   & 2/4        & 3/7     &       &         & 3/7  & 1/1    & 0/2   &         &          & 0/1    & 3/4  & 3/7     &       & 4/7     &          \\
Sampras    & 20/34  & 12/19  &      & 2/2     & 16/20   &          & 8/14   & 0/1     &         & 4/9    & 11/13      & 2/3     & 5/8   & 3/3     & 3/4  & 9/11   &       &         &          & 12/16  & 2/2  & 1/3     & 3/7   &         & 2/3      \\
Wilander   & 2/7    & 3/10   & 0/1  & 5/5     &         &          & 11/20  &         &         &        & 1/2        &         & 7/22  & 6/13    &      & 2/2    &       & 1/1     &          & 2/3    &      &         &       & 1/3     &      \\   
\midrule[1pt]
\multicolumn{26}{l}{\emph{Note:} Due to the space limitation, the names of tennis players have been abbreviated to the first three letters in the horizontal row and their full names can be seen in the vertical column.}\\
\end{tabular*}
\end{sidewaystable}
\restoregeometry

\newgeometry{left=3cm,bottom=2.5cm,top=2cm}
\begin{sidewaystable}
\centering
\caption{Incompete pairwise comparison matrix}
\label{IPCM}
 \tiny
\begin{tabular*}{1\textwidth}{@{\extracolsep{\fill}}llllllllllllllllllllllllll}
\toprule[1pt]
&Aga & Bec & Bor & Con & Cou & Djo & Edb &Fed & Fer & Hew & Kaf & Kue & Len & McE & Moy & Mus & Nad & Nas & New & Raf & Rio & Rod & Saf & Sam & Wil  \\ \hline
Agassi &  & 1.39 &  & 1.07 & 0.90 & &	1.17 & 0.76 & 0.95 & 1.00 & 1.24 & 1.17 & 0.80 & 1.00 & 1.12 & 1.05 & 0.93 &  &  & 1.31 & 0.95 & 1.28 & 1.00 & 0.73 & 1.18   \\
Becker &0.72 &  &  & 1.38 & 1.38 &  & 2.28 &  &  & 1.04 & 1.18 &  & 0.95 & 1.43 & 1.00 & 1.05 &  & 1.03 &  & 1.05 & 1.05 &  & 0.97 & 0.77 & 1.24  \\
Borg &  &  &  & 1.45 &  &  &  &  &  &  &  &  & 1.25 & 1.00 &  &  &  & 1.31 & 0.89 &  &  &  &  &  & 1.03  \\
Connors & 0.93 & 0.73 & 0.69 &  & 0.88 &  & 1.00 &  &  &  &  &  & 0.66 & 0.73 &  &  &  & 0.86 & 1.00 &  &  &  &  & 0.93 & 0.78  \\
Courier & 1.11 & 0.72 &  & 1.13 &  &  & 1.11 &  &  &  & 0.78 & 1.03 & 0.83 & 1.05 & 1.05 & 1.11 &  & 0.97 &  & 0.88 & 0.88 &	 & 1.00 & 0.49 &  \\
Djokovic &  &  &  &  &  &  & & 0.95 & 1.06 & 1.38 &  &  &  &  & 1.00 &  & 0.77 &  &  &  &  & 0.95 & 0.93 &  &   \\
Edberg & 0.85 & 0.44 &  & 1.00 & 0.90 &  &  &  &  &  & 0.95 &  & 1.05 & 0.95 & 1.03 & 1.89 &  &  & & 1.13 & 1.03 & &  & 0.90 & 0.90  \\
Federer & 1.32 &  &  &  &  & 1.05 &  &  & 1.49 & 1.72 & 0.90 & 0.95 &  & & 1.48 & & 0.52 &  &  & 0.88 & 1.07 & 3.31 & 1.64 & 1.03  \\
Ferrero &1.05 &  &  &  &  & 0.95 &  & 0.67 &  & 0.90 & 0.95 & 1.05 &  &  & 1.11 &  & 0.75 &  &  & 1.05 & 1.12 & 0.78 & 1.00 & 0 &   \\
Hewitt & 1.00 & 0.97 &  &  & & 0.72 &  & 0.58 & 1.11 &  & 1.49 & 1.12 &  &  & 1.11 &  & 0.90 &  &  & 1.12 & 1.05 & 1.00 & 1.00 & 1.05 &   \\
Kafelnikov & 0.81 & 0.90 &  &  & 1.28 &  & 1.05 & 1.11  & 1.05 & 0.67 &  & 0.90 &  &  & 1.00 & 0.84 & &  &  & 1.05 & 1.25 & & 1.00 & 0.57 & 1.00  \\
Kuerten & 0.85 &  &  &  & 0.97 &  &  & 1.05 & 0.95 & 0.89 & 1.11 &  &  &  & 1.05 & 1.13 &  &  &  & 1.00 & 1.00 & 1.00 & 1.05 & 0.95 &   \\
Lendl & 1.25 & 1.05 & 0.80 & 1.52 & 1.20 &  & 0.95 &  &  &  &  &  &  & 1.36 & & 1.19 &  & 1.03 &  & 0.97 &  &  &  & 0.90 & 1.54  \\
McEnroe & 1.00 & 0.70 & 1.00 & 1.36 & 0.95 &  & 1.05 &  &  &  &  &  & 0.73 &  &  &  &  & 1.17 & 1.00 &  &  &  &  & 0.88 & 1.05   \\
Moya & 0.89 & 1.00 &  &  & 0.95 & 1.00 & 0.97 & 0.67 & 0.90 & 0.90 & 1.00 & 0.95 &  &  &  & 1.00 & 0.80 &  &  & 1.12 & 0.85 & 0.84 & 1.05 & 0.89 &   \\
Muster & 0.95 & 0.95 &  &  & 0.90 &  & 0.53 &  & &  & 1.19 & 0.88 & 0.84 &  & 1.00 &  &  &  &  & 0.88 & 1.12 &  & 0.97 & 0.65 & 0.93  \\
Nadal & 1.07 &  &  &  &  & 1.29 &  & 1.91 & 1.34 & 1.11 &  &  &  &  & 1.25 &  & &  &  &  & & 1.24 & 1.07 &  &   \\
Nastase &  & 0.97	 & 0.77 & 1.17 & 1.03 &  & &  &  &  &  &  & 0.96 & 0.85 &  &  &  &  & 1.19 &  &  &  &  &  & 0.97  \\
Newcombe &  &  & 1.12 & 1.00 & & & & & &  & & & & 1.00 & & & & 0.84 &  & & &  &  &  &  \\
Rafter & 0.77 & 0.95 &  &  & 1.13 &  & 0.88 & 1.13 & 0.95 & 0.89 & 0.95 & 1.00 & 1.04 &  & 0.89 & 1.13 &  &  &  &  & 1.05 &  & 1.03 & 0.64 & 0.95 \\
Rios &1.05 & 0.95 &  &  & 1.13 &  & 0.97 & 0.93 & 0.89 & 0.95 & 0.80 & 1.00 &  &  & 1.18 & 0.89 &  &  &  & 0.95 &  & 0.93 & 0.89  & 0.93 &   \\
Roddick & 0.78 & &  &  &  & 1.05 &  & 0.30 & 1.28 & 1.00 &  & 1.00 & & & 1.19 &  & 0.80 &  &  &  & 1.07 &  & 1.05 & 1.05 &  \\
Safin & 1.00 & 1.03 &  &  & 1.00 & 1.07 &  & 0.61 & 1.00 & 1.00 & 1.00 & 0.95 &  &  & 0.95 & 1.03 & 0.93 & &  & 0.97 & 1.12 & 0.95 &  & 1.05 &    \\
Sampras & 1.36 & 1.30 & & 1.07 & 2.04 &  & 1.11 & 0.97 &  & 0.95 & 1.77 & 1.05  & 1.11 & 1.13 & 1.12 & 1.53 &  &  &  & 1.57 & 1.07 & 0.95 & 0.95 &  & 1.05  \\
Wilander & 0.85 & 0.80 & 0.97 & 1.28 & & & 1.11 & & & & 1.00 & & 0.65 & 0.95 & & 1.07 && 1.03 & & 1.05  &  &  &  & 0.95 &   \\
\midrule[1pt]
\multicolumn{26}{l}{\emph{Note:} Due to the space limitation, the names of tennis players have been abbreviated to the first three letters in the horizontal row and their full names can be seen in the vertical column.}\\
\end{tabular*}
\end{sidewaystable}
\restoregeometry

\newgeometry{left=3cm,bottom=2.2cm,top=1.8cm}
\begin{sidewaystable}
\centering
\caption{Complete pairwise comparison matrix}
\label{PCM}
 \tiny
\begin{tabular*}{1\textwidth}{@{\extracolsep{\fill}}llllllllllllllllllllllllll}
\toprule[1pt]
    & Aga  & Bec  & Bor  & Con  & Cou  & Djo  & Edb  & Fed  & Fer  & Hew  & Kaf  & Kue  & Len  & McE  & Moy  & Mus  & Nad  & Nas  & New  & Raf  & Rio  & Rod  & Saf  & Sam  & Wil  \\ \hline
Agassi & 1.00 & 1.39 & 0.89 & 1.07 & 0.90 & 1.00 & 1.17 & 0.76 & 0.95 & 1.00 & 1.24 & 1.17 & 0.80 & 1.00 & 1.12 & 1.05 & 0.93 & 1.06 & 1.08 & 1.31 & 0.95 & 1.28 & 1.00 & 0.73 & 1.18 \\
Becker & 0.72 & 1.00 & 0.99 & 1.38 & 1.38 & 0.96 & 2.28 & 0.83 & 1.04 & 1.04 & 1.18 & 1.04 & 0.95 & 1.43 & 1.00 & 1.05 & 0.81 & 1.03 & 1.26 & 1.05 & 1.05 & 1.04 & 0.97 & 0.77 & 1.24 \\
Borg & 1.13 & 1.01 & 1.00 & 1.45 & 1.25 & 1.08 & 1.22 & 0.94 & 1.15 & 1.06 & 1.19 & 1.16 & 1.25 & 1.00 & 1.20 & 1.30 & 0.92 & 1.31 & 0.89 & 1.19 & 1.17 & 1.19 & 1.13 & 0.99 & 1.03 \\
Connors & 0.93 & 0.73 & 0.69 & 1.00 & 0.88 & 0.84 & 1.00 & 0.72 & 0.89 & 0.82 & 0.94 & 0.90 & 0.66 & 0.73 & 0.92 & 1.00 & 0.72 & 0.86 & 1.00 & 0.92 & 0.89 & 0.93 & 0.86 & 0.93 & 0.78 \\
Courier & 1.11 & 0.72 & 0.80 & 1.13 & 1.00 & 0.89 & 1.11 & 0.74 & 0.92 & 0.83 & 0.78 & 1.03 & 0.83 & 1.05 & 1.05 & 1.11 & 0.75 & 0.97 & 1.04 & 0.88 & 0.88 & 0.92 & 1.00 & 0.49 & 0.94 \\
Djokovic & 1.00 & 1.04 & 0.92 & 1.20 & 1.12 & 1.00 & 1.10 & 0.95 & 1.05 & 1.38 & 1.11 & 1.04 & 0.95 & 1.11 & 1.00 & 1.12 & 0.77 & 1.09 & 1.14 & 1.09 & 1.08 & 0.95 & 0.93 & 0.93 & 1.10 \\
Edberg & 0.85 & 0.44 & 0.82 & 1.00 & 0.90 & 0.91 & 1.00 & 0.79 & 0.96 & 0.86 & 0.95 & 0.97 & 1.05 & 0.95 & 1.03 & 1.89 & 0.75 & 0.91 & 0.97 & 1.13 & 1.03 & 0.97 & 0.94 & 0.90 & 0.90 \\
Federer & 1.32 & 1.21 & 1.06 & 1.39 & 1.36 & 1.05 & 1.26 & 1.00 & 1.49 & 1.72 & 0.90 & 0.95 & 1.10 & 1.29 & 1.48 & 1.32 & 0.52 & 1.27 & 1.32 & 0.88 & 1.07 & 3.31 & 1.64 & 1.03 & 1.23 \\
Ferrero & 1.05 & 0.96 & 0.87 & 1.12 & 1.09 & 0.95 & 1.04 & 0.67 & 1.00 & 0.90 & 0.95 & 1.05 & 0.90 & 1.03 & 1.11 & 1.07 & 0.75 & 1.03 & 1.07 & 1.05 & 1.12 & 0.78 & 1.00 & 0.82 & 1.03 \\
Hewitt & 1.00 & 0.97 & 0.95 & 1.22 & 1.20 & 0.72 & 1.17 & 0.58 & 1.11 & 1.00 & 1.49 & 1.12 & 0.97 & 1.14 & 1.11 & 1.17 & 0.90 & 1.10 & 1.17 & 1.12 & 1.05 & 1.00 & 1.00 & 1.05 & 1.13 \\
Kafelnikov & 0.81 & 0.90 & 0.84 & 1.07 & 1.28 & 0.90 & 1.05 & 1.11 & 1.05 & 0.67 & 1.00 & 0.90 & 0.83 & 0.97 & 1.00 & 0.84 & 0.74 & 0.98 & 1.02 & 1.05 & 1.25 & 1.00 & 1.00 & 0.57 & 1.00 \\
Kuerten & 0.85 & 0.96 & 0.86 & 1.11 & 0.97 & 0.96 & 1.03 & 1.05 & 0.95 & 0.89 & 1.11 & 1.00 & 0.89 & 1.03 & 1.05 & 1.13 & 0.78 & 1.02 & 1.06 & 1.00 & 1.00 & 1.00 & 1.05 & 0.95 & 1.03 \\
Lendl & 1.25 & 1.05 & 0.80 & 1.52 & 1.20 & 1.05 & 0.95 & 0.91 & 1.12 & 1.03 & 1.20 & 1.12 & 1.00 & 1.36 & 1.15 & 1.19 & 0.91 & 1.03 & 1.18 & 0.97 & 1.13 & 1.15 & 1.10 & 0.90 & 1.54 \\
McEnroe & 1.00 & 0.70 & 1.00 & 1.36 & 0.95 & 0.90 & 1.05 & 0.78 & 0.97 & 0.88 & 1.03 & 0.97 & 0.73 & 1.00 & 1.00 & 1.10 & 0.78 & 1.17 & 1.00 & 1.01 & 0.96 & 1.00 & 0.93 & 0.88 & 1.05 \\
Moya & 0.89 & 1.00 & 0.84 & 1.08 & 0.95 & 1.00 & 0.97 & 0.67 & 0.90 & 0.90 & 1.00 & 0.95 & 0.87 & 1.00 & 1.00 & 1.00 & 0.80 & 0.98 & 1.03 & 1.12 & 0.85 & 0.84 & 1.05 & 0.89 & 0.99 \\
Muster & 0.95 & 0.95 & 0.77 & 1.00 & 0.90 & 0.89 & 0.53 & 0.76 & 0.93 & 0.86 & 1.19 & 0.88 & 0.84 & 0.91 & 1.00 & 1.00 & 0.74 & 0.91 & 0.94 & 0.88 & 1.12 & 0.93 & 0.97 & 0.65 & 0.93 \\
Nadal & 1.07 & 1.23 & 1.08 & 1.39 & 1.34 & 1.29 & 1.33 & 1.91 & 1.34 & 1.11 & 1.36 & 1.29 & 1.10 & 1.28 & 1.25 & 1.35 & 1.00 & 1.28 & 1.33 & 1.34 & 1.32 & 1.24 & 1.07 & 1.12 & 1.32 \\
Nastase & 0.94 & 0.97 & 0.77 & 1.17 & 1.03 & 0.92 & 1.10 & 0.79 & 0.97 & 0.91 & 1.02 & 0.98 & 0.96 & 0.85 & 1.02 & 1.10 & 0.78 & 1.00 & 1.19 & 1.01 & 0.99 & 0.99 & 0.96 & 0.80 & 0.97 \\
Newcombe & 0.93 & 0.79 & 1.12 & 1.00 & 0.96 & 0.88 & 1.03 & 0.76 & 0.94 & 0.86 & 0.98 & 0.94 & 0.85 & 1.00 & 0.97 & 1.06 & 0.75 & 0.84 & 1.00 & 0.97 & 0.94 & 0.96 & 0.91 & 0.82 & 0.96 \\
Rafter & 0.77 & 0.95 & 0.84 & 1.09 & 1.13 & 0.92 & 0.88 & 1.13 & 0.95 & 0.89 & 0.95 & 1.00 & 1.04 & 0.99 & 0.89 & 1.13 & 0.74 & 0.99 & 1.03 & 1.00 & 1.05 & 1.02 & 1.03 & 0.64 & 0.95 \\
Rios & 1.05 & 0.95 & 0.85 & 1.12 & 1.13 & 0.92 & 0.97 & 0.93 & 0.89 & 0.95 & 0.80 & 1.00 & 0.89 & 1.04 & 1.18 & 0.89 & 0.76 & 1.01 & 1.07 & 0.95 & 1.00 & 0.93 & 0.89 & 0.93 & 1.00 \\
Roddick & 0.78 & 0.96 & 0.84 & 1.08 & 1.08 & 1.05 & 1.03 & 0.30 & 1.28 & 1.00 & 1.00 & 1.00 & 0.87 & 1.00 & 1.19 & 1.07 & 0.80 & 1.01 & 1.04 & 0.98 & 1.07 & 1.00 & 1.05 & 1.05 & 1.00 \\
Safin & 1.00 & 1.03 & 0.88 & 1.16 & 1.00 & 1.07 & 1.06 & 0.61 & 1.00 & 1.00 & 1.00 & 0.95 & 0.91 & 1.07 & 0.95 & 1.03 & 0.93 & 1.04 & 1.10 & 0.97 & 1.12 & 0.95 & 1.00 & 1.05 & 1.05 \\
Sampras & 1.36 & 1.30 & 1.01 & 1.07 & 2.04 & 1.07 & 1.11 & 0.97 & 1.22 & 0.95 & 1.77 & 1.05 & 1.11 & 1.13 & 1.12 & 1.53 & 0.89 & 1.25 & 1.22 & 1.57 & 1.07 & 0.95 & 0.95 & 1.00 & 1.05 \\
Wilander & 0.85 & 0.80 & 0.97 & 1.28 & 1.06 & 0.91 & 1.11 & 0.81 & 0.97 & 0.88 & 1.00 & 0.97 & 0.65 & 0.95 & 1.01 & 1.07 & 0.76 & 1.03 & 1.05 & 1.05 & 1.00 & 1.00 & 0.96 & 0.95 & 1.00 \\
\midrule[1pt]
\multicolumn{26}{l}{\emph{Note:} Due to the space limitation, the names of tennis players have been abbreviated to the first three letters in the horizontal row and their full names can be seen in the vertical column.}\\
\end{tabular*}
\end{sidewaystable}
\restoregeometry


Firstly,
the statistical competition results matrix (i.e. with additive preference relations) is converted into the pairwise comparison matrix (i.e. with multiplicative preference relations) following the method described in \cite{bozoki2016application}. 
Fig.\ref{transformation} illustrates this transformation process. The converted pairwise comparison matrix is shown in Tab.\ref{IPCM}.

Secondly, DEMATEL-based completion method is adopted to complete the pairwise comparison matrix. Following the procedure of the proposed method, the complete pairwise comparison matrix is obtained, which is shown in Tab.\ref{PCM}.

Thirdly, based on the complete pairwise comparison matrix, the priorities of tennis players can be calculated via (\ref{weights of AHP}). 
Tab.\ref{Ranks} shows the final ranking of the 25 tennis players.

\begin{table}[htbp]
\centering
\caption{Ranking of 25 tennis players}
\label{Ranks}
\scriptsize
\begin{tabular*}{1\textwidth}{@{\extracolsep{\fill}}clcclc}
\toprule[1pt]
\textbf{Rank} & Player   & Priority \quad & \textbf{Rank}  & Player     & Priority \\ \hline
\textbf{1}    & Nadal    & 0.0503 & \textbf{14}   & Rios       & 0.0380 \\
\textbf{2}    & Federer  & 0.0503 & \textbf{15}   & McEnroe    & 0.0380 \\
\textbf{3}    & Sampras  & 0.0467 & \textbf{16}   & Nastase    & 0.0380 \\
\textbf{4}    & Borg     & 0.0444 & \textbf{17}   & Rafter     & 0.0378 \\
\textbf{5}    & Lendl    & 0.0436 & \textbf{18}   & Wilander   & 0.0378 \\
\textbf{6}    & Becker   & 0.0430 & \textbf{19}   & Kafelnikov & 0.0375 \\
\textbf{7}    & Hewitt   & 0.0414 & \textbf{20}   & Edberg     & 0.0374 \\
\textbf{8}    & Djokovic & 0.0412 & \textbf{21}   & Moya       & 0.0371 \\
\textbf{9}    & Agassi   & 0.0409 & \textbf{22}   & Newcombe   & 0.0365 \\
\textbf{10}   & Safin    & 0.0392 & \textbf{23}   & Courier    & 0.0360 \\
\textbf{11}	  & Kuerten  & 0.0390 & \textbf{24}   & Muster     & 0.0353 \\
\textbf{12}   & Roddick  & 0.0384 & \textbf{25}   & Connors    & 0.0339 \\
\textbf{13}   & Ferrero  & 0.0383 &      &            &    	\\ 
\bottomrule[1pt] 
\end{tabular*}
\end{table}

\section{Discussion}

The result in last section indicates that the proposed method has ability to complete the pairwise comparison matrices in AHP.  
In this section, we would further analyze the superiority and efficiency of the proposed method.

The proposed method has the following advantages:
\begin{itemize}
\item
The proposed method estimates the missing values in pairwise comparison matrix from a new perspective (i.e. by calculating the total relations among alternatives).
\item
The proposed method has low computational cost.
\end{itemize}

After reviewing the prior studies, an assessment for the computational cost of the existing completion methods and the proposed method is done. The methods for comparison satisfy:

\begin{itemize}
\item
The methods belong to completion methods, which have ability to estimate the missing values in pairwise comparison matrix.
\item
The pairwise comparison matrices, which can be completed by these completion methods, present reciprocal preference relations among alternatives. In other words, the matrices request that the multiplication of each pair of symmetric values along the diagonal line is equal to 1.
\end{itemize}

As for these completion methods, some of them estimate the missing values through iteration, or based on the complex theories/methods, which causes the methods expensive computational cost.
For example, 
Bozo\'iki \emph{et al.}\cite{bozoki2010optimal} developed a completion method which estimates the missing values in pairwise comparison matrix through iteration and optimizes the matrix consistency  based on graph theory. 
Gomez \emph{et al.}\cite{gomez2010estimation} proposed a completion method based on the neural network. Neural network as one of the branches of artificial intelligence leads to the complexity of the completion method.
As for the completion methods without iteration or complex theories/algorithms such as the methods listed in Tab.\ref{comparative methods}, time complexity is used to accurately reflect their computational cost. Assume there is an $m \times m$ incomplete pairwise comparison matrix, which has n missing values, the time complexity of these methods is shown in Tab.\ref{timeComplexity}.

\begin{table}[htbp]
\centering
\caption{Methods for comparison}
\label{comparative methods}
\scriptsize
\newcommand{\tabincell}[2]{\begin{tabular}{@{}#1@{}}#2\end{tabular}}
\begin{tabular*}{1\textwidth}{@{\extracolsep{\fill}}lll}
\toprule[1pt]
&  	\quad \quad \quad \quad \quad \quad \quad \quad \quad \quad Description &  \quad \quad Source 	\\ \hline
Method 1    &	\tabincell{l}{
A completion method extends the geometric mean induced bias matrix to est-\\~~imate the missing judgments and improves the consistency ratios by the least\\~~absolute error method and the least square method.
}
&  Ergu \emph{et al.}\cite{ergu2016estimating} \\
Method 2 	&  	\tabincell{l}{A consistency-based completion procedure, which takes a measure to maintain \\ ~~the consistency level of matrices before and after completion.} 	& Alonso \emph{et al.}\cite{alonso2008consistency}  \\
Method 3 	&  \tabincell{l}{A completion method to obtain the most consistent/similar matrix with the in-\\~~itial incomplete matrix based on Frobenius norm.} &  Ben\'itez \emph{et al.}\cite{benitez2015consistent} \\
\bottomrule[1pt] 
\end{tabular*}
\end{table}

\begin{table}[htbp]
\centering
\caption{Time complexity of methods}
\label{timeComplexity}
\footnotesize
\begin{tabular*}{1\textwidth}{@{\extracolsep{\fill}}lcccc}
\toprule[1pt]
					& Method 1 			& Method 2  			& Method 3 		& The proposed method \\ 
\hline
Time complexity		& {$m^3$}			& {$m^3$}		& {$m^2+mn$}			& {$m^2$}				\\  
\bottomrule[1pt]
\end{tabular*}
\end{table}

It can be seen from the Tab.\ref{timeComplexity} that the time complexity of the proposed method is lower than the others.
In addition, the most existing completion methods have the following two operations: 
\begin{itemize}
\item Divide the matrix into upper and lower triangular matrix. 
The pairwise comparison matrix in AHP has a property that the multiplication of each pair of symmetric values along the diagonal line is equal to 1. These methods firstly focus on the completion of upper/lower triangle of incomplete matrix and use the property to accomplish the other triangular matrix.
\item Traverse the pairwise comparison matrix to obtain the positions of unavailable values, then estimate these missing values in sequence. 
\end{itemize}

However, the proposed method completes the matrix from a holistic viewpoint without dividing the matrix into the upper and lower triangular matrix.
Furthermore, based on DEMATEL, all the missing values in matrix can be simultaneously estimated.
Hence, following the procedure of the proposed method, the completion of incomplete pairwise comparison matrix becomes straightforward and efficient.


\section{Conclusion}
\label{conclusion}
Pairwise comparison matrix plays a pivotal role in AHP. However, in many cases, only partial information in pairwise comparison matrix is available, which obstructs the subsequent operations of the classical AHP.
In this paper, we propose a new completion method for incomplete pairwise comparison matrix in AHP. 
Based on DEMATEL, the proposed method has ability to complete the matrix by obtaining the total (direct and indirect) relations from the direct preference relations between each pair of alternatives. 
Furthermore, an application of the proposed method for ranking tennis players is presented. 
From the ranking result and the comparison with several representative completion methods,
the proposed method is proved to be effective and efficient. 
The proposed method provides a new perspective to complete the matrix with explicit physical meaning. Besides, the proposed method has low computation cost. 
This promising method has a wide application in multi-criteria decision-making. In our further study, we would extend the proposed method for the incomplete fuzzy pairwise comparison matrices.

\section*{Acknowledgment}

The work is partially supported by National Natural Science Foundation of China (Grant Nos. 61573290,61503237), China State Key Laboratory of Virtual Reality Technology and Systems, Beihang University (Grant No.BUAA-VR-14KF-02).

\section*{References}
%


\bibliographystyle{elsarticle-num}
\bibliography{references}









\end{document}